\documentclass[11pt,leqno,svgnames]{amsart}
\usepackage{amssymb}
\usepackage{mathrsfs}
\usepackage{color}
\usepackage{ulem}
\usepackage{xcolor}
\usepackage{dsfont}
\usepackage{fancyhdr}
\numberwithin{equation}{section}

\newcommand\R{{\mathbb R}}

\def\AA{{\mathcal A}}
\def\BB{{\mathcal B}}

\def\DD{{\mathcal D}}

\def\II{{\mathcal I}}
\def\JJ{{\mathcal J}}
\def\KK{{\mathcal K}}
\def\LL{{\mathcal L}}
\def\MM{{\mathcal M}}

\def\OO{{\mathcal O}}
\def\PP{{\mathcal P}}
\def\QQ{{\mathcal Q}}
\def\RR{{\mathcal R}}

\def\TT{{\mathcal T}}
\def\UU{{\mathcal U}}
\def\VV{{\mathcal V}}
\def\WW{{\mathcal W}}
\def\XX{{\mathcal X}}
\def\YY{{\mathcal Y}}
\def\ZZ{{\mathcal Z}}

\def\BBB{{\mathscr{B}}}

\def\Ppp{{\mathbf P}}

\newtheorem{theo}{Theorem}[section]

\newtheorem{lem}[theo]{Lemma}
\newtheorem{cor}[theo]{Corollary}

\newcommand{\beqn}{\begin{equation}}
\newcommand{\eeqn}{\end{equation}}
\newcommand{\bear}{\begin{eqnarray}}
\newcommand{\eear}{\end{eqnarray}}
\newcommand{\bean}{\begin{eqnarray*}}
\newcommand{\eean}{\end{eqnarray*}}
\newcommand{\bal}{\begin{aligned}}
\newcommand{\eal}{\end{aligned}}

\newcommand{\be}{\begin{equation}}
\newcommand{\ee}{\end{equation}}
\newcommand{\ba}{\begin{aligned}}
\newcommand{\ea}{\end{aligned}}

\newcommand{\Black}{\color{black}}

\def\signsm{\bigskip \begin{center} {\sc St\'ephane Mischler\par\vspace{3mm}
Universit\'e Paris-Dauphine, PSL Research University,
 \par
CNRS, UMR [7534], CEREMADE,
\par
Place du Mar\'echal de Lattre de Tassigny
75775 Paris Cedex 16\par
FRANCE\par
e-mail:} \tt{mischler@ceremade.dauphine.fr} \end{center}}

\def\signcq{\bigskip \begin{center} {\sc Cristobal Qui{\~n}inao\par\vspace{3mm}
 Universidad de O'Higgins,  \par 
 Intituto de Ciencias de la Ingenieria,  \par 
 Avenida Libertador Bernardo O'Higgins 611,  
 \par 
 Rancagua,  CHILE\par

e-mail:} \tt{cristobal.quininao@uoh.cl} \end{center}}

\def\signqw{\bigskip \begin{center} {\sc Qilong Weng\par\vspace{3mm}
Universit\'e Paris-Dauphine, PSL Research University,
 \par
CNRS, UMR [7534], CEREMADE,
\par
Place du Mar\'echal de Lattre de Tassigny
75775 Paris Cedex 16\par
FRANCE\par
e-mail:} \tt{weng@ceremade.dauphine.fr} \end{center}}

\begin{document}

\title{Weak { and strong} connectivity regimes for \\ a  general time elapsed neuron network model}

%

\author{S. Mischler, C. Qui\~ninao, Q. Weng}

\begin{abstract}
For large fully connected neuron networks, we study the dynamics of
homogenous assemblies of interacting neurons described by time
elapsed models.

Under general assumptions on the  firing rate which include the ones made
in previous works \cite{PK1,PK2,MSWQ}, we establish accurate estimate on  the long time behavior
of the solutions in the weak  and the strong connectivity regime both in the case with and without delay.
Our results improve  \cite{PK1,PK2} where a less accurate estimate was established
and \cite{MSWQ} where only smooth firing rates were considered.

Our approach combines several arguments introduced in the above previous works
as well as a slightly refined version of the Weyl's and spectral mapping theorems
presented in \cite{Voigt80,MS}.

%
\end{abstract}

\maketitle

\begin{center} {\bf Version of \today}
\end{center}

\vspace{0.3cm}

\bigskip
\bigskip


\tableofcontents



\maketitle

\section{Introduction}
\label{sec:Intro}

%
%
%
%

The information transmission and processing mechanism in the nervous
systems relies on the quantity of electrical pulses as the reflect to
incoming stimulations, during which the neurons experience a period of
recalcitrance called discharge time before reactive  (for more information about neuronal networks and mean-field approach see e.g.~\cite{MR2902610,MR2974499}). In this work, we
shall focus on the model describing the neuronal dynamics in
accordance with this kind of discharge time which has been introduced
and studied in \cite{gerstner2002spiking,PK1,PK2}. In order to show the response to the
recovery of the neuronal membranes after each discharge, the model
consider an instantaneous firing rate depending on the time elapsed
since last discharge as well as the inputs of neurons. This sort of
models are also regarded as a mean field limit of finite number of
neuron network models referred to \cite{MR3311484,FL*,RT*,MR3449317,MR3718099,MR3631217,Q*}.

For a local time (or internal clock) $x\ge0$ corresponding to the
elapsed time since the last discharge, we consider the dynamic of the
neuronal network with the density number of neurons $f=f(t,x)\geq0$
in state $x \ge 0$ at time $t \ge 0$, given by the following nonlinear
time elapsed (or of age structured type) evolution equation
\begin{subequations}\label{eq:ASM}
 \begin{align}
   &\partial_t f=-\partial_x f-k(x,\lambda \,  m(t))f=:\mathcal{L}_{\lambda m(t)}f,\\
   &f(t,0)=p(t), \ \ f(0,x)=f_0(x),
 \end{align}
\end{subequations}
where $k(x,\lambda \, \mu)\ge0$ denotes the firing rate of a neuron in
state $x$ and in an environment $\mu\ge0$ formed by the global
neuronal activity with a network connectivity parameter $\lambda\ge0$
corresponding to the strength of the interactions. The total density
of neurons $p(t)$ undergoing a discharge at time $t$ is defined through
\begin{subequations}\label{eq:disch}
 \begin{gather}
   p(t):=\mathcal{P}[f(t); m(t)],\\
   \intertext{where}
   \PP [g,\mu] = \PP_\lambda[g,\mu] :=  \int_0^\infty k(x,\lambda \mu) g(x)\mathrm{d}x,\label{eq:foP}
 \end{gather}
\end{subequations}
 while the global neuronal activity $m(t)$ at time $t\geq0$ taking into
account the interactions among the neurons resulting from earlier
discharges is given by
\beqn\label{eq:delay}
m(t):=\int_0^\infty p(t-y)b(\mathrm{d}y).
\eeqn
Here the delay distribution $b$ is a probability measure taking into account
the persistence of the electric activity to those discharges in the
network. In the sequel, we will consider the two following situations
respectively:

$\bullet$ The {\it case without delay} when $b=\delta_0$, and then $m(t) = p(t)$;

$\bullet$ The {\it case with delay} when $b$ is a smooth function.

\medskip
Observe that the solution $f$ of the time elapsed equation
\eqref{eq:ASM} satisfies
$$
\frac{\mathrm{d}}{\mathrm{d}t}\int_0^{\infty}f(t,x)\mathrm{d}x=f(t,0)-
\int_0^{\infty}k(x,\lambda m(t))f(t,x)\mathrm{d}x=0,
$$
in both  cases. That implies the conservation of the total density
number of neurons 
which can be thus normalized to $1$. As a consequence, we assume in the sequel
\beqn\label{eq:MassConservation}
\langle f(t,\cdot) \rangle = \langle f_0 \rangle = 1,
\quad \forall t\geq0, \quad \langle g \rangle := \int_0^{\infty} g(x)\mathrm{d}x.
\eeqn

\smallskip
We call steady state a couple $(F_\lambda,M_\lambda)$ of a nonnegative function and a positive real number which satisfies
\begin{subequations}\label{eq:StSt}
 \begin{align} \label{eq:StSt1}
   &0=-\partial_x F_\lambda - k(x,\lambda \,  M_\lambda) F_\lambda= \mathcal{L}_{\lambda M_\lambda} F_\lambda,\\
   &F_\lambda(0)= M_\lambda, \quad \langle F_\lambda \rangle = 1.
 \end{align}
\end{subequations}
Noticing that the associated network activity and the discharge
activity are equal constants for a steady state  because $\langle b\rangle=1$.

\smallskip
Our main purpose in this paper is to prove
existence, uniqueness and exponential asymptotic stability of solutions to the time elapsed
evolution equation \eqref{eq:ASM} in weak and strong connectivity regimes,
which is a range of connectivity parameter $\lambda\in (0,\lambda_0)\cup(\lambda_{_\infty},\infty)$, with $\lambda_0 > 0$ small enough  and $\lambda_{_\infty}>0$ large enough, chosen in such a way that the
nonlinear term in equations \eqref{eq:ASM} and \eqref{eq:StSt}  is not too strong.

 \smallskip
 These results are obtained for a rather large class of firing rate. More precisely, we make the physically reasonable assumptions
\beqn\label{hyp:a0}
k \ge0, \quad\partial_x k\ge0, \quad k'=\partial_\mu k\geq0,
\eeqn
meaning that neurons increases their firing rate as the global activity is higher and as the elapsed time since last decharge increases, 
\beqn\label{hyp:a1}
0<k_0:=\lim_{x\rightarrow\infty}k(x,0)\leq\lim_{x,\,\mu\rightarrow\infty}k(x,\mu)=:k_{_1}<\infty,
\eeqn
as well as the regularity assumption
\beqn\label{hyp:a2}
K(x,\cdot):=\int_0^x k(y,\cdot)\mathrm{d}y\in C^0(\R_+), \quad\forall x>0.
\eeqn
We will also need the stronger regularity assumption $k \in   \hbox{\rm Lip}_\mu L^1_x $.  For the weak connectivity regime, we assume that for some
 $\xi>0$ small enough and for any $\mu_0 >0$, there exists $\lambda_0>0$ small enough such that 
\beqn\label{hyp:a3}
\int_0^\infty |k(x,\lambda\mu_2) - k(x,\lambda\mu_{_1})| \, dx \le \xi \, |\mu_2 - \mu_{_1}|,\quad\forall\mu_{_1},\mu_2\in(0,\mu_0), \,\, \forall \, \lambda\in(0,\lambda_0).
\eeqn
While in the strong connectivity regime, we assume that for some the same $\xi > 0$ as in the assumption \eqref{hyp:a3} and for any $\mu_{_\infty}>0$, there exists $\lambda_{_\infty}>0$ large enough such that
\beqn\label{hyp:a4}
\int_0^\infty\big|k(x,\lambda\,\mu_2)-k(x,\lambda\,\mu_{_1})\big|\,\mathrm{d}x\le\xi\,|\mu_2 - \mu_{_1}|,\quad\forall\mu_{_1},\mu_2\in(\mu_{_\infty},\infty),
\,\,\lambda\in(\lambda_{_\infty},\infty).
\eeqn
A possible example of firing rate which fulfills the above condition \eqref{hyp:a3} for the weak connectivity regime  is the ``step function firing rate"  considered in \cite{PK1,PK2} which is given by
\beqn\label{hyp:s1}
k(x,\mu) = {\bf 1}_{x > \sigma(\mu)},\quad \sigma' \le 0,
\eeqn
\beqn\label{hyp:s2}
\sigma_+ := \sigma(0), \quad\sigma_-:= \sigma(\infty), \quad\sigma_-<\sigma_+<1,
\eeqn
where $\sigma$ satisfies the regularity condition
\beqn\label{hyp:s3}
\quad \sigma,\sigma^{-1}\in W^{1,\infty}(\R_+).
\eeqn
Similarly,  the above condition \eqref{hyp:a4} for the strong connectivity regime is met for a ``step function firing rate"
introduced in  \cite{PK1,PK2} given by the same function as above which additionnaly fullfils
$$
s \, |\sigma'(s)| \to 0 \quad\hbox{as}\quad s \to\infty.
$$
%
In the case with delay, we assume that the delay distribution is associated to a measurable function, namely
$b(\mathrm{d}y) = b(y)\mathrm{d}y$ with $b \in L^1(\R_+)$, and satisfies the exponential bound
\begin{equation}\label{hyp:del}
  \exists\delta>0, \quad\int_0^\infty e^{\delta
y}\, b(y) 
\, \mathrm{d}y <\infty.
\end{equation}

Our first result establishes the existence and  uniqueness of weak solution to the evolution problem \eqref{eq:ASM}. We call weak solution a function $0 \le f\in C(\R_+;L^1(\R_+)_w)\cap L^\infty(\R_+^2)$ such that
$$
\int_0^T\!\!\!\int_0^\infty f \, (\partial_t \varphi + \partial_x \varphi) \, dxdt = \int_0^T\!\!\! \int_0^\infty k(x,\lambda m) f \varphi \, dxdt -  \int_0^T p \, \varphi(t,0) \, dt
$$
for any $\varphi \in C^1_c(\R_+^2)$, where $p$ and $m$ satisfy  \eqref{eq:disch}-\eqref{eq:delay}. Here and below $L^1_q(\R_+)$, for $q > 0$, stands for the space of $L^1$ functions $f$ such that $x^q f \in L^1$ and $L^1_w$ denotes  the $L^1(\R_+)$ space
endowed with the weak topology $\sigma(L^1,L^\infty)$.

\begin{theo}\label{th:EaU}
We consider a firing rate $k$ satisfying \eqref{hyp:a0}-\eqref{hyp:a1}-\eqref{hyp:a2} and a initial datum $0 \le f_0\in L^1_{ q}(\R_+)\cap L^\infty(\R_+)$, for some $q > 0$, with 
total density number of neuron $1$.
We further assume that one of the following conditions holds:

(1) the delay distribution  $b$  satisfies \eqref{hyp:del};

(2) $b = \delta_0$, $k$ satisfies  \eqref{hyp:a3}  and $\lambda \in (0,\lambda_0)$, for $\lambda_0 > 0$ small enough;

(3) $b = \delta_0$, $k$ satisfies \eqref{hyp:a4}  and $\lambda \in (\lambda_{_\infty},\infty)$, for $\lambda_{_\infty} > 0$ large enough, as well as
\beqn\label{eq:kf0}
 \kappa_0:=\int_0^\infty k(x,0)f_0(x)\,\mathrm{d}x>0.
\eeqn
%
In any of these three cases,
there exists a 
weak solution $0 \le f\in C(\R_+;L^1(\R_+)_w)\cap L^\infty(\R_+;L^1_q(\R_+))\cap L^\infty(\R_+^2)$ to the evolution equation \eqref{eq:ASM}-\eqref{eq:disch}-\eqref{eq:delay} for some functions $m,\,p\in C(\R_+)$ which satisfies the total number density of neurons conservation \eqref{eq:MassConservation} as well as 
\beqn\label{eq:th:L1infty}
\| f_t \|_{L^\infty} \le \| f _0 \|_{L^\infty}   +
k_{_1},\quad \|  f\|_{L^1_q }\le \| f_0 \|_{L^1_q } + K_q, 
\quad \forall \, t \ge 0, \eeqn
for some constant $K_q = K_q(k) \ge 0$ and
\beqn\label{eq:mLinfty}
0 \le \kappa_1 
\le\|m_t \|_{L^\infty}\le k_{_1}, \quad \forall \, t \ge \tau, 
\eeqn
with $ \kappa_1 > 0$ when $\kappa_0 > 0$ or when $\tau > 0$ is large enough. The solution is furthermore unique in case (2) and (3).
\end{theo}


 Our proof is based on a Schauder fixed point theorem in the case with delay  (1) and on a Banach fixed point theorem in cases without delay (2) and (3). 

\medskip
As a second step, we state an existence of solution to the stationary problem  \eqref{eq:StSt} and the uniqueness
of that solution in the weak and strong connectivity regime.

\begin{theo}\label{th:SS}
Under the above assumption \eqref{hyp:a0}-\eqref{hyp:a1}-\eqref{hyp:a2} on the firing rate,  for any $\lambda\ge0$, there exists at least one couple $(F_\lambda(x),M_\lambda)\in W^{1,\infty}(\R_+)\times\R_+$ solution to the stationary problem \eqref{eq:StSt}, and  such that
\beqn\label{ineq:EoF}
0\le F_\lambda(x)\lesssim e^{{-k_0\over2}x},\quad |F'_\lambda(x)|\lesssim e^{{-k_0\over2}x},\quad x\ge0.
\eeqn
Moreover, when we assume additionaly that  \eqref{hyp:a3} and \eqref{hyp:a4} hold, there exist $\lambda_0>0$ small enough  and $\lambda_{_\infty}$ large enough, such that the above steady state is unique for any $\lambda\in[0,\lambda_0)\cup(\lambda_{_\infty},\infty]$.
\end{theo}

The proof being identical to the ones presented in \cite{PK1} and \cite[Theorem~2.1]{MSWQ}, it will be skipped.

\medskip
Finally our third and main result in the present paper states the exponential nonlinear stability of the above stationary state in the weak  and strong connectivity regime.

\begin{theo}\label{th:MR} We assume that $k$, $b$ and $f_0$ satisfy the same conditions (1), (2) or (3) as in Theorem~\ref{th:EaU} and furthermore $k$ satisfies  \eqref{hyp:a3} and \eqref{hyp:a4}. There exist $\lambda_0>0$ small enough, $\lambda_\infty> 0$ large enough, some constants $\alpha< 0$ and $C \ge 1$ such that for any  $\lambda\in(0, \lambda_0) \cup(\lambda_\infty,+\infty)$ the solution  $f$ to the evolution equation  \eqref{eq:ASM}-\eqref{eq:disch}-\eqref{eq:delay} built in Theorem~\ref{th:EaU} furthermore satisfies
\beqn\label{th:NLStab}
\|f(t,.)-F_\lambda \|_{L^1}\leq C  e^{\alpha t}, \qquad \forall \, t \ge 0 .
\eeqn
\end{theo}

\smallskip
This theorem generalizes  to the delay case  the similar results obtained in \cite{PK1,PK2}  and it
generalizes the similar result obtained in \cite{MSWQ} to a more general firing rate including the step function rate considered in \cite{PK1,PK2}.

The proof is mainly based on an extension of the abstract semigroup theory developed in \cite{MS,MSWQ} which has probably its own interest. It uses an auxiliary linear problem introduced in
 \cite{PK1,PK2} instead of the linearized equation  considered in \cite{MSWQ}. Both arguments together make possible to get ride of the smoothness assumption needed in  \cite{MSWQ} and moreover, allow us to consider the large connectivity regime, and also generalize the stability results established in \cite{PK1,PK2}.


\smallskip
Our approach is thus quite different from the usual way to deal with delay equations which consists in using the specific framework of   ``fading memory space", which goes back at least to Coleman \& Mizel   \cite{MR0210343},  or the theory of ``abstract algebraic-delay differential systems" developed by O. Diekmann and co-authors \cite{MR1345150}. It is also different from the previous works \cite{PK1,PK2} where the asymptotic stability analysis were performed by taking advantage of the ``step function firing rate" \eqref{hyp:s1}, making possible to find a suitable norm such that the problem becomes dissipative. 

\smallskip

\medskip
This paper is organized as follows.
In Section~\ref{sec:ExitSteS}, we establish the existence and uniqueness results for the evolution equation as stated in Theorem~\ref{th:EaU}.
The estimate on the long time behavior of solutions as formulated in Theorem~\ref{th:MR} is established in  Section~\ref{sec:WithoutDelay} in the case without delay. The case with delay is tackled in Section~\ref{sec:WithDelay}.

%
%

\section{Existence of solutions}\label{sec:ExitSteS}

\subsection{Delay case}

In order to establish the existence of a solution to \eqref{eq:ASM}-\eqref{eq:disch}-\eqref{eq:delay}, we will apply a Schauder fixed point argument. To begin with, we analyze the continuity property of the functional $\PP$ defined in \eqref{eq:foP}.

\begin{lem}\label{lem:cop}
Assume \eqref{hyp:a0}-\eqref{hyp:a1}-\eqref{hyp:a2}.
Consider a sequence $(m_n)$ of nonnegative real numbers converging to a limit $m$ in $\R$
as well as a sequence of  functions $(f_n)$ which converges  to $f$ in the sense of the weak topology $\sigma(L^1,L^\infty)$ and is uniformly bounded in $L^\infty$. We then have
$$
\PP[f_n,m_n]\to\PP[f,m], \quad as \ n\to\infty.
$$
\end{lem}
\proof {\sl Step 1. Continuity of $k$.} We are going to show that
\beqn\label{eq:amkTOam}
k(\cdot,m_n)\to k(\cdot,m), \quad a.e. \quad \hbox{as} \ n\to\infty.
\eeqn
We first assume that $(m_n)$ is increasing. The sequence $(k(\cdot,m_n))$ is also increasing because of assumption \eqref{hyp:a0}. Moreover, since $k$ is bounded from assumption \eqref{hyp:a1},  there exists some $\bar{k}(x)$ such that
$$
k(x,m_n)\to\bar{k}(x), \quad \hbox{as} \ n\to\infty, \quad \hbox{for any} \, x\ge0,
$$
which in turn implies
$$
K(x,m_n)=\int_0^x k(y,m_n)\mathrm{d}y\to\int_0^x\bar{k}(y)\mathrm{d}y.
$$
From assumption \eqref{hyp:a2}, we deduce
$$
\int_0^x k(y,m)\mathrm{d}y=K(x,m)=\int_0^x\bar{k}(y)\mathrm{d}y, \quad\forall x\ge0.
$$
Thus, we clearly have
$$
\bar{k}(x) = k(x,m), \quad \hbox{for a.e.} \ x \ge 0.
$$
The same holds in the case when $(m_n)$ is a decreasing sequence. In the general case,  we may define two monotonous sequences $(m^i_n)$, $i = 1,2$, such that
$m^1_n \le m_n \le m^2_n$ for any $n \ge 1$ and such that $m^i_n\to m$ as $n\to\infty$ for $i=1,2$. Then $k(x,m^1_n) \le k(x,m_n) \le k(x,m^2_n)$ for any $n \ge 1$, $x \ge 0$ and
 $k(x,m^i_n) \to k(x,m)$ as $k \to \infty$ for a.e. $x \ge 0$ and for $i=1,2$. We immediately conclude that \eqref{eq:amkTOam} holds.

\smallskip
\noindent{\sl Step 2. Continuity of the functional $\PP$.} We compute
\bean
\PP[f_n,m_n]-\PP[f,m] &=& \int_0^\infty k(\cdot,m_n) f_n-\int_0^\infty k(\cdot,m) f\\
 &=& \int_0^\infty(k(\cdot,m_n)-k(\cdot,m))f_n+\int_0^\infty k(\cdot,m)(f_n-f)\\
 &:=& I_{_1}+I_2.
\eean
\Black From the assumption \eqref{hyp:a0} and the weak
convergence of $f_n$ in $L^1$, we have $I_2\to0$, as $n\to\infty$.
We write \bean I_{_1} =
\int_0^R(k(x,m_n)-k(x,m))f_n(x)\mathrm{d}x+\int_R^\infty(k(x,m_n)-k(x,m))f_n(x)\mathrm{d}x.
\eean From Step 1 and the assumption that $(f_n)$ is bounded in
$L^\infty$ and uniformly integrable at the infinity (as a
consequence of its weak $\sigma(L^1,L^\infty)$ convergence and the
Dunford-Pettis theorem), we deduce that \bean |I_{_1}| \le
\|f_n\|_{L^\infty}\int_0^R|k(x,m_n)-k(x,m)|\mathrm{d}x+2k_{_1}\int_R^\infty
f_n(x)\mathrm{d}x\to0, \eean as $R\to\infty$  and  $n\to\infty$. The
two above estimates togeter imply the conclusion. \endproof


\smallskip
In a next step, we fix $T > 0$ and we analyse the linear mapping which associates to a given function $m\in C([0,T])$ the solution $f\in C([0,T];L^1)\cap L^\infty([0,T];L^\infty)$ to the transport equation
\beqn\label{eq:Pb2}
 \begin{aligned}
 & \partial_t f+\partial_x f+k(x,m(t))f=0\\
 & f(t,0)=\PP[f,m(t)], \quad f(0,x)=f_0(x).
 \end{aligned}
\eeqn
The following lemma gives the continuity of this mapping.

\begin{lem}\label{lem:cof}
Fix $T>0$. Consider  a sequence $(m_n)$ such that $m_n\to \bar m$ in $C([0,T])$, as $n\to\infty$. There exists then a sequence $(f_n)$ of solutions to the linear transport equation \eqref{eq:Pb2} associated to $(m_n)$ and this one satisfies
$$
f_n\to \bar f\quad in \ C([0,T];L^1_w)\cap
L^\infty([0,T];L^\infty{ \cap L^1_q}), \quad \hbox{as} \
n\to\infty,
$$
where $\bar f$ stands for the solution to the linear transport equation \eqref{eq:Pb2} associated to $\bar m$. \end{lem}


\proof {\sl Step 1. Existence of $f_n$.} For any $m\in C([0,T])$ and
any $0 \le g_{_1} \in X_T := C([0,T];L^1) \cap
L^\infty(0,T;L^\infty{ \cap L^1_q})$ we may associate $0 \le g_2
\in X_T$ as the solution to the equation \beqn\label{eq:Pb1}
 \begin{aligned}
 & \partial_t g_2+\partial_x g_2+k(x,m(t))g_2=0\\
 & g_2(t,0)=p_{_1}(t):=\PP[g_{_1},m(t)], \quad g_2(0,x)=f_0(x),
 \end{aligned}
\eeqn which is classically defined through the characteristic
method. More precisely, we introduce the space
$$
\mathscr{C}:=\{0\le g\in X_T; \,  \|g\|_{L^1}\le   e^{k_{_1}t},
\|g\|_{L^\infty}\le\|f_0\|_{L^\infty}+k_{_1},  \|g\|_{L^1_q}\le
C_q\|f_0\|_{L^1_q}+e^{k_1 T}\},
$$
for some  constant $C_q \ge 1$ that we specify below  and we consider $g_{_1} \in \mathscr{C}$.   Integrating the equation \eqref{eq:Pb1} on $x$, we find
$$
\frac{\mathrm{d}}{\mathrm{d}t}\int_0^\infty g_2(t,x)\mathrm{d}x=\int_0^\infty(g_{_1}(t,x)-g_2(t,x))k(x,m(t))\mathrm{d}x\le k_{_1}\|g_{_1}\|_{L^1}\le k_{_1}e^{k_{_1}t},
$$
which implies
$$
\|g_2\|_{L^1}\le\int_0^t k_{_1}e^{k_{_1}s}\mathrm{d}s+\|f_0\|_{L^1}=e^{k_{_1}t}.
$$
 We  define 
$$
\omega(x):=  {\bf 1}_{x \le x_0} + {x^q \over x_0^q} \, {\bf
1}_{x > x_0} , 
$$
with $x_0 \ge 1$ large enough such that $q/x- k(\mu,x) \le q/x_0 - k(0,x_0) \le - k_0/2 \le 0$ for any $x \ge x_0$ and $\mu \ge 0$.
Similarly as above, we deduce 
\bean
\frac{\mathrm{d}}{\mathrm{d}t}\int_0^\infty g_2\omega &=&
\int_0^\infty k(x,m(t)) g_1+\int_0^\infty g_2 [ \partial_x\omega -  k(x,m(t)) \omega] \\
&\le& k_1e^{k_1t} , 
\eean
which implies
$$
\|g_2\|_{L^1_\omega }
\le e^{k_1t} - 1 + \|f_0\|_{L^1_\omega} \le e^{k_1t} + x_0^{-q}  \|f_0\|_{L^1_q}. 
$$
We finally get 
$$
\|g_2\|_{L^1_q }
\le  C_q e^{k_1t} +    \|f_0\|_{L^1_q},
$$
with $C_q := x_0^q$. 
 Thanks to the method of characteristics the solution
$g_2(t,x)$ to equation \eqref{eq:Pb1} can be expressed as
$$
g_2(t,x)=\left\{
\begin{aligned}
 f_0(x-t) e^{-\int_0^t
  k(s+x-t,m(s))\mathrm{d}s}, & \quad & \forall\, x\geq t,\\
 p_{_1}(t-x) e^{-\int_0^x
  k(s,m(s+t-x))\mathrm{d}s}, & \quad & \forall\, x\leq t,
\end{aligned}
\right.
$$
which implies
$$
\|g_2\|_{L^\infty}\le\|f_0\|_{L^\infty}+\|p_{_1}\|_{L^\infty}\le\|f_0\|_{L^\infty}+k_{_1}.
$$
Denoting $g_2 :=\II (g_{_1})$, we have proved $\II(\mathscr{C})\subset\mathscr{C}$. On the other hand, denoting $h_{_1} = g_{_1} -\tilde g_{_1}$ and $h_2 := \II(g_{_1}) - \II(\tilde g_{_1})$ for $g_{_1}, \tilde g_{_1} \in \mathscr{C}$, a similar computation as above leads to
$$
\sup_{t \in [0,T]} \| h_2(t) \|_{L^1} \le {1 \over k_{_1}} \Bigl( e^{k_{_1} T} - 1 \Bigr) \sup_{t \in [0,T]} \| h_{_1}(t) \|_{L^1},
$$
from what we conclude to the existence of a unique function $ f
\in\mathscr{C}$ such that $\II (f) = f$ by a classical {contraction}
fixed point Theorem for $T > 0$ small enough. We get $T > 0$
arbitrary by iterating the argument. We thus get the existence of
the sequence $(f_n)$ and the possible limit $\bar f$ by applying the
above construction with $m = m_n$ and $m = \bar m$.

 \smallskip\noindent{\sl Step 2. Uniform estimates on
$\bar{f}$.} By integrating the transport equation \eqref{eq:Pb2}, we
obtain that the solution concerves the total density number of neurons as given by \eqref{eq:MassConservation}. For any solution ${\bar f}$ to
the equation \eqref{eq:Pb2}, we deduce
\bean
\frac{\mathrm{d}}{\mathrm{d}t}\int_0^\infty {\bar f}\omega 
&=&
\int_0^\infty k(x,m(t)) {\bar f} +\int_0^\infty {\bar f} [ \partial_x\omega -  k(x,m(t)) \omega] \\
&\le& k_1 - {k_0\over2} \int_{x_0}^\infty {\bar f}   \omega  \le k_1 + {k_0\over2}  - {k_0\over2}  \int_0^\infty {\bar f} \omega,
\eean
which implies
$$
\|\bar f\|_{L^1_\omega }\le \min \bigl\{ 1 + \frac{2k_1}{k_0}, \| f_0 \|_{L^1_\omega } \bigr\}. 
$$
We finally define 
$$
\bar{\mathscr{C}}:=\{0\le f\in
X_T;\|f\|_{L^1}=1,\|f\|_{L^\infty}\le\|f_0\|_{L^\infty}+k_1,\|f\|_{L^1_q}\le\|f_0\|_{L^1_q}+ K_q\},
$$
with $K_q := 2 x_0^q ( 1 + k_1/k_0)$. 
By construction,   $\bar{f}$ obviously belongs to
$\bar{\mathscr{C}}$.

\smallskip\noindent{\sl Step 3. Continuity of the mapping.} From  equation \eqref{eq:Pb2}, we observe that
$$
 \partial_t f_n\le-\partial_x f_n, \quad f_n (t,0)\le k_{_1}, \quad f {_n}(0,x)=f_0,
$$
which implies
$$
f_n(t,x)\le k_{_1}\mathbf{1}_{x\le t}+f_0(x-t)\mathbf{1}_{x\ge t},
$$
for any $t\le T$.  From this upper bound, one gets
$$
\int_R^\infty f_n\le\int_R^\infty f_0(x-t)\mathrm{d}x,
$$
for any $R\ge T$, and in particular
$$
\int_R^\infty f_n\le \int_{R-T}^\infty f_0\to0,
$$
as $R\to\infty$ and uniformly in $n \ge 1$. From equation \eqref{eq:Pb1}, for any $\varphi \in C^1_c( (0,\infty))$, we also have
$$
{d \over dt} \int_0^\infty f_n \varphi \, dx = A^\varphi_n,
$$
with $A^\varphi_n$ bounded in $L^\infty(0,T)$ uniformly in $n \ge 1$.
Together with the fact that $f_n \in \mathscr{C}$ for any $n \ge 1$, we may use the de la Vallee-Poussin Theorem and the Dunford-Pettis Lemma to conclude that there exists $f\in C([0,T];L^1_w)\cap L^\infty([0,T];L^\infty)$ and a subsequence $f_{n'}$ of the sequence $f_n$ such that $f_{n'}\rightharpoonup f$ weakly.
We deduce  $\PP[f_{n'},m_{n'}]\to\PP[f,m]$ as ${n'}\to\infty$, from Lemma \ref{lem:cop}. We finally conclude by passing ${n'}$ to the limit in the equation \eqref{eq:Pb2} with $m_{n'}$.
\endproof

\proof[Proof of Theorem \ref{th:EaU} - the delay case] We recall that $b \in L^1(\R_+)$ in that case.  We consider the application $\JJ : C([0,T])\to C([0,T])$,
defined as
$$
\JJ(m)(t):=\int_0^t p(t-y)b(\mathrm{d}y), \quad \forall \, m \in C([0,T]), \,\, \forall \, t \in [0,T],
$$
where  $p(t)=\PP[f(t,\cdot),m(t)]$ and $f \in\bar{\mathscr{C}}$
is a solution to \eqref{eq:Pb2} which existence has been established
during the proof of Lemma~\ref{lem:cof}. From Lemma \ref{lem:cop},
we deduce that the application $m\to\PP[f,m]$ is continuous and so
is $\JJ$. Define
$$
\mathscr{K}:=\{m\in C([0,T]), \|m\|_{L^\infty}\le k_{_1}\}.
$$
Obviously, $\mathscr{K}$ is a convex subset of $C([0,T])$ and, for any $m\in\mathscr{K}$, we have
$$
\|\JJ(m)\|_{L^\infty}\le\int_0^t|b(y)|\|p\|_{L^\infty}\mathrm{d}y\le k_{_1},
$$
so that $\JJ : \mathscr{K} \to \mathscr{K}$. On the other hand, for any $\epsilon>0$, there exists $\theta>0$ such that for any $t,s\in[0,T]$ satisfying $|t-s|<\theta$, we have
\bean
|\JJ(m)(t)-\JJ(m)(s)| &\le& \int_0^\infty|b(t-y)-b(s-y)|\|p\|_{L^\infty}\mathrm{d}y\\
 &\le& \|b(t-\cdot)-b(s-\cdot)\|_{L^1(\R_+)}k_{_1}\|f_0\|_{L^1}\\
 &\le& k_{_1}\|\tau_\theta b-b\|_{L^1(\R_+)}<\epsilon,
\eean where $\tau_\theta b:=b(\cdot+\theta)$, which implies that
$\JJ$ is equicontinuous. Thanks to the Arzela-Ascoli Theorem, we
deduce that $\JJ(\mathscr{K})$ is compactly embedded into
$\mathscr{K}$. Using the Schauder-Brouwer fixed point Theorem, the
application $\JJ$ admits a fixed point $m\in\mathscr{K}$. The
corresponding solution $f$ to the equation \eqref{eq:Pb2} is also a
solution to the equation
\eqref{eq:ASM}-\eqref{eq:disch}-\eqref{eq:delay}. Iterating on $T$,
we deduce the existence of a global solution $(f,m,p)$ to equation
\eqref{eq:ASM}-\eqref{eq:delay}, which satisfies the bound in
Theorem~\ref{th:EaU}. 

 In order to prove the lower bound \eqref{eq:mLinfty} on $m$, we recall that $x_0 \ge 1$ has been defined such that 
$$
k(x,\mu)\ge \frac{k_0}2\mathbf{1}_{x\ge x_0}, \quad \forall \, x > 0, \, \mu >0.
$$
For $t\le x_0$, thanks to the characteristics expression, we have
\bean
p(t) &\ge& \int_0^\infty k(x,0)f_0(x-t) e^{-\int_0^t k(s+x-t,m(s))\,\mathrm{d}s}\,\mathrm{d}x\\
 &\ge& \int_0^\infty k(x,0)f_0(x)e^{-k_1 t}\,\mathrm{d}x\\
 &\ge&   e^{-k_1x_0}\kappa_0.
\eean
We consider now the case $t> x_0$. Directly from \eqref{eq:ASM}, we have 
$$
\partial_t f\ge\partial_x f-k_1 f, 
$$
which implies 
$$
f(t,x)\ge e^{-k_1x_0}f(t-x_0,x-x_0).
$$
We then deduce
\bean
p(t) &\ge& e^{-k_1x_0}\int_0^\infty \frac{k_0}2\mathbf{1}_{x\ge x_0}f(t-x_0,x-x_0)\,\mathrm{d}x\\
 &\ge& \frac{k_0}2e^{-k_1x_0}\int_0^\infty f(t-x_0,x)\,\mathrm{d}x\\
 &=&\frac{k_0}2e^{-k_1x_0}.
\eean
All together, we deduce the same lower bound for $m(t)$ from the definition \eqref{eq:delay}.
 \Black
\endproof


\subsection{Without delay case}

We will need the following auxilliary result.
We define the function $\Phi :  L^1(\R_+) \times \R \to \R$ by
$$
\Phi [g,\mu] := \int_0^\infty k(x,\lambda \mu) g (x) \, \mathrm{d}x - \mu.
$$

\begin{lem}\label{lem:varphig} Assume \eqref{hyp:a0}-\eqref{hyp:a1}-\eqref{hyp:a2}-\eqref{hyp:a3}-\eqref{hyp:a4}. For any bounded set $\mathscr{C} \subset {\Ppp} \cap L^\infty$, there exists $\lambda_0 = \lambda_0(\mathscr{C}) > 0$ and $\lambda_\infty=\lambda_\infty(\mathscr{C})>0$, such that for any $\mu_0>0$ and $\mu_\infty>0$, there exists a function $\varphi_\lambda : \mathscr{C} \to \R$ which is Lipschitz continuous  in the sense of the $L^1$ norm and such that $\mu = \varphi_\lambda[g]$ is the unique solution to the equation
$$\ 
\begin{aligned}
& \mu \in(0,\mu_0), & \quad \Phi(g,\mu) = 0, &\quad\forall\lambda\in(0,\lambda_0),\\
& \mu \in(\mu_\infty,k_{_1}), & \quad \Phi(g,\mu) = 0, &\quad\forall\lambda\in(\lambda_\infty,\infty).
\end{aligned}
$$
\end{lem}

\proof[Proof of Lemma \ref{lem:varphig}] The proof is similar to the one of \cite[Lemma~2.8]{MSWQ}, thus we skip the existence part and we present the uniqueness part here. Fix $R > 0$ and take $f,g \in L^1 \cap L^\infty$ and $\mu,\nu \in  (0,\mu_0)$ or $\mu,\nu\in(\mu_\infty,k_{_1})$ such that
$$
\|f \|_{L^\infty} \le R \quad \hbox{ and } \quad \Phi(f,\mu) =  \Phi(g,\nu) = 0.
 $$
We have
$$
\nu - \mu = \int_0^\infty k(x,\lambda\, \nu) (g - f) \, dx +   \int_0^\infty \big(k(x,\lambda\, \nu) -  k(x,\lambda\, \mu) \big) f  \, dx  ,
$$
with
$$
\Bigl| \int_0^\infty k(x,\lambda\, \nu) (g - f) \, dx  \Bigr| \le k_{_1} \, \|f-g\|_{L^1}.
$$
 From the assumption \eqref{hyp:a3}-\eqref{hyp:a4} and the uniform estimate on $f$, there holds
$$
\Bigl| \int_0^\infty \big(k(x,\lambda\,\nu) -  k(x,\lambda\,\mu)\big) f   \, dx \Bigr| \le R \, \xi\,|\mu-\nu|,\quad\forall\lambda\in[0,\lambda_0)\cup(\lambda_{_\infty},\infty].
$$
We then deduce
$$
|\mu - \nu|\,\le2k_{_1}\|f-g\|_{L^1},
$$
for any $\lambda \in [0,\lambda_0)\cup(\lambda_{_\infty},\infty]$, with $\lambda_0=\lambda_0(R) > 0$ small enough and $\lambda_{_\infty}=\lambda_{_\infty}(R)>0$ large enough. 
That implies the uniqueness of the solution $\mu  = \varphi_\lambda(f)\in \R$ to the constraint problem $\Phi(f,\mu) = 0$ for any given $f \in L^1 \cap L^\infty$ and the Lipschitz continuity of $\varphi_\lambda$.

\proof[Proof of Theorem \ref{th:EaU} - The  case without delay.] We fix $\lambda_0, \lambda_{_\infty} > 0$ as defined in Lemma~\ref{lem:varphig}. For a given function $m \in C([0,T])$, we define $M (t) := \varphi_\lambda(f(t,.)) \in C([0,T])$, where $f \in C([0,T];L^1(\R_+)) \cap L^\infty([0,T] \times\R_+)$ is the solution of \eqref{eq:Pb2} associated to $m$. We denote $\II(m) := M$. For two given $m_{_1},m_2 \in C([0,T))$, we denote $f_{_1}, f_2 \in C([0,T];L^1(\R_+)) \cap L^\infty([0,T] \times\R_+)$ the associated solutions to \eqref{eq:Pb2} and we easily compute
\bean
{\mathrm{d} \over \mathrm{d}t} \int |f_2 - f_{_1}|
&\le& 2 \int |k(m_2) f_2 - k(m_{_1}) f_{_1}|
\\
&\le& 2 \| f_0 \|_{L^\infty} \xi\, |m_2 - m_{_1}| + 2 k_{_1} \int |  f_2 -   f_{_1}|,\quad\forall\lambda\in(0,\lambda_0)\cup(\lambda_{_\infty},\infty).
\eean
We deduce that  $m \mapsto f$ is Lipschitz from $C([0,T])$ to $C([0,T];L^1(\R_+))$ with constant $CT$. As a consequence,  $\II$ is Lipschitz from $C([0,T])$ into itself with constant $C'T$.  Choosing $T > 0$ small enough, the mapping $\II$ is a contraction and admits a unique fixed point thanks to the Banach fixed point theorem.  Iterating on $T$, we deduce the existence and uniqueness of a global solution $(f,m)$ to equation \eqref{eq:ASM}-\eqref{eq:disch} in the case without delay in both weak  and strong connectivity regimes.
\endproof

\section{A Weyl's and spectral mapping theorm}
\label{sec:W&SMtheorems}

In this section we establish a simple version of Weyl's and spectral mapping theorem for semigroup
in an abstract setting which slightly generalizes the versions of the same theorems
established in \cite{MS,MSWQ}. More precisely, we consider the generators $\LL$ and $\BB$
of two semigroups $S_\LL$ and $S_\BB$ in a Banach space $\XX$. We denote $\AA := \LL - \BB$  as well as $R_\LL(z) := (\LL-z)^{-1}$
$R_\BB(z) := (\BB-z)^{-1}$ the resolvent operators defined in the corresponding resolvent sets. We  assume
that for some fixed $a^* \in \R$ the following growth and regularizing estimates hold true for any $a > a^*$:
  \begin{itemize}
   \item[{\bf (H1)}] $\BB$ is $\AA$-power dissipative in $X$, in the sense that
\beqn\label{eq:hypH1}
\forall \, \ell \ge 0, \quad t \mapsto \| S_\BB *  (\AA S_\BB)^{(*\ell)}(t ) \|_{\BBB(\XX)} \, e^{-a t} \in L^\infty (0,\infty)
\eeqn
and $u := S_\BB \AA$ satisfies
\beqn\label{eq:hypH1bis}
\exists \, n \ge 1, \,\, C \in (0,\infty), \,\, \forall \, f \in \XX, \qquad \int_0^\infty \| u^{(*n)}(t) f \|_\XX \, e^{-at} \, dt \le C\, \| f \|_\XX.
\eeqn
\item[{\bf (H2)}]  For the same integer $n \ge 1$, the operator $\UU := - R_\BB\AA $ is  power regular in the sense that
\beqn\label{eq:hypH2}
\exists \, \alpha > 0, \,\,  C \in (0,\infty), \,\,   \qquad   \| \UU(z)^n  \|_{\BBB(\XX)}  \le C \, \langle x \rangle^{-\alpha}, \quad \forall \, z \in \Delta_a
\eeqn
and
\beqn\label{eq:hypH2bis}
\forall \, M > 0, \,\, \exists \,  C \in (0,\infty), \,\,   \qquad   \| \UU(z)^n  \|_{\BBB(\XX,\YY)}  \le C , \quad \forall \, z \in \Delta_a \cap B(0,M),
\eeqn
for some linear space $\YY$ such that the embedding $\YY \subset \XX$ is compact.
\end{itemize}

For a given operator $L$ we denote $\Sigma(L)$ its spectral set and we define $\Sigma_d(L)$ the discrete spectrum as the set of isolated eigenvalues with finite dimensional associated eigenspace. We also denote $D(L)$ the domain and $RL$ the range.

 \begin{theo}  \label{theo:GalWeylSG}
We make the above growth and regularizing assumptions {\bf (H1)} and  {\bf (H2)} on $\AA$ and $\BB$ for some $a^* < 0$ and we assume furthermore that
$$
\Sigma (\LL) \cap \bar\Delta_0 = \{ 0 \} \subset \Sigma_d(\LL),
$$
so that there exist a finite rank projector $\Pi_0 \in \BBB(\XX)$ and an operator $T_0 \in \BBB(R \Pi_0)$ satisfying
$\LL \Pi_0 = \Pi_0 \LL = T_0 \Pi_0$, $\Sigma(T_0) = \{ 0 \}$.
There exist  $a < 0$ and $C \ge 1$  such that
  \begin{equation}\label{bddSlambda}
     \forall \, t \ge 0, \quad
     \Bigl\| e^{t \, \LL} - e^{t \, T_0} \, \Pi_0
     \Bigr\|_{\mathscr{B}(\XX)}
     \le C  \, e^{a \, t}.
   \end{equation}

\end{theo}

%
For two given time dependent operators valued functions $U$ and $V$, we define the convolution product
$$
(U*V)(t) := \int_0^t U(t-s) V(s) \, ds.
$$
We also denote $V^{(*1)} = V$ and $V^{(*\ell)} := V * V^{(* (\ell-1))}$ for any $\ell \ge 2$.

\smallskip\noindent
{\sl Proof of Theorem~\ref{theo:GalWeylSG}. Step 1.  }  We define
$$
\VV(z) := \RR_\BB(z) - \dots + (-1)^{n-1} \RR_\BB(z) \, (\AA \RR_\BB(z))^{n-1}
$$
and
$$\WW(z) := (-1)^{n} (  \RR_\BB(z)\AA)^{n},
$$
where  $n \ge 1$ is  the integer given by assumption {\bf (H1)}. From the definition $\LL = \AA + \BB$, we immediately have
$$
R_\LL = R_\BB - R_\BB \AA R_\LL,
$$
and by iterating that relation, we deduce
$$
R_\LL =  \VV + \WW \, R_\LL,
$$
or equivalently
$$
 (I - \WW) R_\LL = \VV.
$$
Thanks to \eqref{eq:hypH2},  for $M$ large enough, we have
$$
z \in \Delta_a, \, |z| \ge M \quad\Rightarrow \quad \| \WW(z) \|_{\BBB(X)} \le {1 \over 2}.
$$
We get that $R_\LL(z) =  ( I - \WW(z))^{-1} \VV(z) $ is well defined and uniformly bounded in the region $\Delta_a \backslash B(0,M)$, or in other word
$\Sigma(\LL) \cap \Delta_a \subset B(0,M)$.

\smallskip
On the other hand, $\Phi := I - \WW$ is holomorphic on $\Delta_{a^*}$ and $R(\WW) \subset \YY \subset\subset \XX$ because of \eqref{eq:hypH1}-\eqref{eq:hypH2bis}.
Together with $\Phi(M)$ is invertible, we may use Ribari{\v{c}}-Vidav-Voigt theory \cite{RVidav,Voigt80} and deduce that $R_\LL =  (I-\WW)^{-1} \VV$ is a degenerate-meromorphic operator and next that   $\Sigma(\LL) \cap \Delta_{a^*} $ is discrete. All together, we have proved that there exists $a < 0$ such that $\bar\Delta_a \cap \Sigma(\LL) = \{0 \}$.

\smallskip
\noindent {\sl Step 2. }    For any integer $N \ge 1$ and  iterating the Duhamel formula
$$
S_\LL\ =   S_\BB   +  ( S_\BB \AA) * S_\LL,
$$
we have \bean
S_\LL \,  (I-\Pi_0)
=   \sum_{\ell=0}^{N-1}  S_\BB *  (\AA S_\BB)^{(*\ell)}   (I-\Pi_0) +   (S_\BB \AA )^{(*N)} * (S_\LL (I-\Pi_0) ).
\eean
For $b >  \Lambda(\LL)$,  we may use the inverse Laplace formula
\bean
\TT(t) f
&:=&  (\AA S_\BB)^{(*N)} * (S_\LL (I-\Pi_0) ) (t)  f
\\
&=& \lim_{M' \to \infty}  {i \over
  2\pi}\int_{b-iM'}^{b+iM'} e^{zt} \,  (-1)^{N+1} \,( \RR_\BB(z)\AA)^{N}   (I-\Pi_0)   \RR_\LL(z) \,   f \, dz,
\eean
for any $f \in D(\LL)$  and $t \ge 0$, and we emphasize that the term $\TT(t) f$ might be only defined as a semi-convergent  integral.
Because $z \mapsto ( \RR_\BB(z)\AA)^{N}   (I-\Pi_0)   \RR_\LL(z) $ is a bounded  analytic function on a neighborhood of $\bar \Delta_a$,
we may move the segment on which the integral is performed, and we obtain
\beqn\label{eq:SL=2}
\TT(t) f =   \lim_{M' \to \infty}  {i \over
  2\Pi_0}\int_{a-iM'}^{a+iM'} e^{zt} \,  (-1)^{N+1} \, ( \RR_\BB(z)\AA)^{N}   \RR_\LL(z)  (I-\Pi_0)    f \, dz,
\eeqn
 for any $f \in D(\LL)$  and $t \ge 0$.   In order to conclude we only have to explain why the RHS term in \eqref{eq:SL=2} is appropriately bounded for $N$ large enough. We define
$$
\WW(z) := \RR_\LL(z) \,  (\AA \RR_\BB(z))^{N}
$$
for $z \in \Delta_a \backslash B(0,M)$,  $N := ( [1/\alpha]+1) n$.
From Step 1 and \eqref{eq:hypH2}, we deduce
\beqn\label{eq:ch5:estimWW}
\| \WW(z)\|_{\BBB(X)}  \le  {C \over | y |^\beta}, \quad \forall \, z = a + y,  \, |y| \ge M,
\eeqn
with  $\beta :=  ( [1/\alpha]+1) \alpha > 1$.
 We then have
\bean
\bigl\| \TT (t) f \bigr\|_{\BBB(X)}
&\le&   {e^{at} \over 2\pi} \int_{a-iM}^{a+iM} \| (R_\BB(z)\AA )^{N} \|_{\BBB(X)} \,  \|  (I-\Pi_0) R_\LL(z) \|_{\BBB(X)} \, dy
\\
&&\quad+ {e^{at} \over 2\pi} \int_{\R \backslash [-M,M]} \|   \WW(a+iy)  (I-\Pi_0)   \|_{\BBB(X)} \, dy,
\eean
where the first integral is finite thanks to  $ \Sigma(\LL  (I-\Pi_0)  ) \cap [a-iM,a+iM] = \emptyset$ and \eqref{eq:hypH2}, while the second
integral is finite because of \eqref{eq:ch5:estimWW}.
\qed

\section{Case without delay}
\label{sec:WithoutDelay}

In this section, we present the proof of our main result Theorem~\ref{th:MR}
in the case without delay.

\subsection{An auxiliary linear equation}
We introduce the  auxiliary linear equation on the variation $g$  given by
\bear
   \nonumber
   &&\partial_t g +\partial_x g+k_\lambda g =0,\\
   &&g(t,0)= \MM_\lambda[g(t,.)],
   \quad g(0,x)=g_0(x),
   \label{eq:ch3:renewal}
\eear
with the notations
\beqn\label{eq:deMMeps}
\MM_\lambda[h]  := \int_0^{\infty}k_\lambda h \, \mathrm{d}x, \quad k_\lambda :=k(x,\lambda M_\lambda),
\eeqn
and where $M_\lambda$ is defined in Theorem~\ref{th:SS}.
The corresponding
linear operator $\LL$ is
$$
\LL g:=-\partial_x g-k_\lambda g
$$
in the domain
$$
D(\LL):=\{g\in W^{1,1}(\R_+), \, g(0)=\MM_\lambda[g]\}
$$
generating the semigroup $S_{\LL}$ in the Lebesgue space $X:=L^1(\nu)$ for some polynomial weight function $\nu :=1 + x^q$, $q > 0$.  For any initial datum $g_0\in X$, the weak solution of the linearized equation is given by $g(t)=S_{\LL}(t)g_0$.
By regarding the boundary condition as a source term, we may rewrite the above equation as
\beqn\label{eq:ASMDL}
\partial_t g=\Lambda  g :=-\partial_x g-k_\lambda g  +\delta_{x=0}\MM_\lambda[g],
\eeqn
with the associated semigroup $S_{\Lambda_\lambda}$, acting on the space of bounded Radon measures
$$
\XX:=M^1(\mathbb{R}_+) = \{ g \in (C_0(\R))'; \,\, \hbox{supp} \, g \subset \R_+ \},
$$
endowed with the weak $*$ topology $\sigma(M^1,C_0)$, where $C_0$
represents the space of continuous functions converging to $0$ at
infinity.

\begin{theo}\label{th:CWR}
For any $\lambda\ge0$, there exist $\alpha<0$ and $C>0$ such that  $\Sigma(\LL ) \cap \Delta_\alpha = \{ 0 \}$ and
\beqn\label{eq:CWR}
\| S_{\LL}(t) g_0 \|_X \leq C\, e^{\alpha t} \, \| g_0 \|_X, \quad \forall \, t \ge 0,
\eeqn
for any $g_0 \in X$, $\langle g_0 \rangle = 0$.
\end{theo}

We proceed in several steps.

\begin{lem}\label{lem:SLge0} The semigroup $S_\LL$ is well defined in $L^1$ and it is positive in the sense that $S_\LL(t) f_0 \ge 0$
for any $f_0 \in L^1$, $f_0 \ge 0$ and any $t \ge 0$.
\end{lem}

The proof being exactly the same as for \cite[Lemma~2.5]{MSWQ} it is skipped.

\begin{lem}\label{lem:LstrongMax}
$-\LL$ satisfies the following version of the strong maximum principle:
 for any given $g\in \XX_+$ and $\mu\in\mathbb{R}$, there holds
$$
g\in D(\LL)\setminus\{0\} \ and \ (-\LL+\mu)g\geq0 \
imply \ g>0.
$$
\end{lem}

The proof being exactly the same as for  \cite[Lemma~2.6]{MSWQ} it is skipped. As an immediate consequence of Theorem~\ref{th:SS},
Lemma~\ref{lem:SLge0} and Lemma~\ref{lem:LstrongMax}, we obtain the following result about the first eigenvalue and eigenspace associated to $\LL$.
We refer to   \cite[Proof of Theorem 2.4.]{MSWQ} or \cite[Proof of Theorem 5.3]{MS,MSerra} where similar results are established.

\begin{cor}\label{cor:FirstEV}
There hold $\Sigma (\LL) \cap \bar\Delta_0 = \{ 0 \} $ and $N(\LL) = \hbox{span} (F_\lambda)$.
\end{cor}

\smallskip
We come to the slightly new argument we need in order to generalize the proof presented in  \cite{MSWQ} to the non smooth firing rate we are considered here. Because of \eqref{hyp:a0}-\eqref{hyp:a1}, we have
$$
k_\lambda \in L^\infty(\R_+), 
\quad k_\lambda (x) \ge k_0/2\, {\bf 1}_{x \ge x_0},
$$
for some $x_0 \in [0,\infty)$, and we set
$$
a^* := - k_0/2 < 0.
$$
 We rewrite the evolution equation as
\beqn\label{eq:ch3:renewallBIS}
\partial_t f = \LL f = \AA f + \BB f,
\eeqn
with $\AA$ and $\BB$ defined by
\bean
(\AA f ) (x) &:=& \delta_{x=0} \KK[f], \quad \KK[f] := \int_0^\infty k_\lambda(y) \, f(y) \, dy
\\
(\BB f ) (x)
&:=& - \partial_x f(x)  - k_\lambda(x) f(x),
\eean
and we emphasize that the boundary condition in  \eqref{eq:ch3:renewal}  has been equivalently replaced by the term $\AA f$ involving a Dirac mass $ \delta_{x=0}$.

\begin{lem}\label{lem:Renewal3} For any $a > a^*$ and any $R > 0$, the operators $\AA$ and $\BB$ satisfy  :

(i) the operator $\BB-a$ is  dissipative in $L^1(\tilde\nu)$, with norm equivalent to the norm of $X$;

(ii) the operators valued function of time $t \mapsto  ( S_\BB \AA)^{(*\ell)} * S_\BB (t)\, e^{-at}$ is bounded in $L^\infty(0,\infty;\BBB(X))$;

(iii) the operators valued holomorphic function $z \mapsto \langle z \rangle (R_\BB (z) \AA)^2(z)$ is  bounded in $\BBB(X)$ uniformly in $z \in \Delta_a$;

(iv) $z \mapsto  (R_\BB (z) \AA)^2(z)$ is  bounded in $\BBB(X,Y)$ uniformly in $z \in \Delta_a \cap B(0,R)$, with $Y := BV \cap L^1_{q+1}$.
\end{lem}

\noindent{\sl Proof of Lemma~\ref{lem:Renewal3}. } During the proof we write $k=k_\lambda$.

\smallskip\noindent {\sl Step 1. }
In order to prove the first point, we fix $a > a^*$ and we introduce the modified weight function
$$
\tilde\nu (x) := {e^{ax} \over e^{ax_{_1}}} \, {\bf 1}_{x \le x_{_1}} +  {x^q \over x_{_1}^q} \, {\bf 1}_{x > x_{_1}},
$$
with $x_{_1} > \max(1,q/(a-a^*))$. We  compute
\bean
\BB^* \tilde\nu &=& \partial_x \tilde\nu - k \tilde\nu \le a \, \tilde\nu \quad\hbox{on}\ [0,x_{_1}]
\\
\BB^* \tilde\nu &\le& \Big( {q \over x} - k \Big) \tilde\nu \le a \, \tilde\nu \quad\hbox{on}\ (x_{_1},\infty),
\eean
from what we deduce
$$
\int_0^\infty (\BB f) f/|f| \tilde\nu  = \int_0^\infty (\BB^* \tilde\nu) |f|  \le a \int_0^\infty |f| \tilde\nu
$$
and $\BB - a$ is dissipative in $L^1(\tilde\nu)$. 

\smallskip\noindent {\sl Step 2. } From the first step we have $\| S_\BB(t) \|_{\BBB(X)} = \OO(e^{at})$.
We deduce (ii) recursively.

\smallskip\noindent {\sl Step 3. }
We have
$$
S_\BB (t) f (x) = f(x-t) \, \exp( K(x-t) - K(x)),
$$
with
$$
K(x) := \int_0^x k(u) \, du.
$$
We deduce successively
\bean
S_\BB (t) \AA f (x)
=  \delta_{x-t=0}  \, \exp( K(x-t) - K(x)) \, \KK[f] ,
\eean
next
\bean
g_{t-s} (x)
&:=& \AA  S_\BB (t-s) \AA f (x)
\\
&=&  \delta_{x=0}  \,  k(t-s)  \exp( - K(t-s)) \, \KK[f]  ,
\eean
and finally
\bean
(S_\BB \AA)^{(*2)}(t) f (x)
&=&  \int_0^t g_{t-s}(x-s) \, \exp( K(x-s)) - K(x)) \, ds
\\
&=& {\bf 1}_{t \ge x}   \,  k(t-x)  \exp( - K(t-x)  - K(x))  \, \KK[f] =: \varphi_t(x) \, \KK[f].
\eean
Summarizing, we have
$$
(S_\BB \AA)^{(*2)}(t) f  =  \varphi_t  \, \KK[f],
$$
with
$$
\varphi_t(x) = \psi(t-x)  \exp(  - K(x)), \quad \psi (u) := {\bf 1}_{u \ge 0} \, k(u) \, e^{-K(u)}.
$$
We then compute
$$
\hat \varphi_z(x) = \hat \psi (z) \, e^{- z x} \, e^{-K(x)}.
$$
On the one hand, using that $K(u) \ge k_0 \, u - \tilde{k}$, for any $a > a^* = - k_0$ and some $\tilde{k} \in \R$, we have
\bean
\| e^{- z x } \, e^{-K(x)} \|_{L^1(\nu)}
&\le& \int_0^\infty  e^{- \Re e z x - K(x)} \, \langle x \rangle^q  \, dx
\\
&\le& \int_0^\infty  C   \, e^{-   (\Re e z -a) x  }   \, dx
\le { C\over  \Re e z-a} , 
\eean
for any $z \in \Delta_a$ and some constant $C \in (0,\infty)$.
On the other hand, when furthermore $k' \in L^1(0,\infty)$, we may perform one integration by parts and we get
\bean
\hat \psi(z)
&=& {1 \over z} \Bigl( k(0) -  \int_0^\infty (k^2(u) - k'(u)) \, e^{- K(u)} \, e^{- zu } \, du \Bigr).
\eean
As a consequence and similarly as above, we have
\bean
|\hat \psi(z) | &\le& {1 \over |z|} \, \Bigl( \| k \|_{L^\infty} +  C_{_1} \| k' \|_{M^1} + C_2 \, {\| k \|_{L^\infty}^2 \over \Re e z-a} \Bigr),
\eean
for any $z \in \Delta_a$, $a > a^*$. By a standard regularization argument, we get the same estimate in the general case when $k' \in M^1([0,\infty)$.
All together, we obtain
\bean
\| (R_\BB (z) \AA)^2(z) f \|_X
&=& \| e^{-z x  - K(x)} \|_X \, |\hat \psi (z) | \, | \KK[f] |
\\
&\le& {C_a \over |z|} \, \| f \|_X ,
\eean
for any  $z \in \Delta_a$, $a > a^*$ and a constant $C_a$ depending of $a$, $L^1(\nu)$ and $k$ (through the quantities $\| k \|_{L^\infty}$, $\| k' \|_{M^1}$ and $k_{_1}$).

\smallskip\noindent {\sl Step 4. }  We observe that
$$
\| \hat\psi_z \|_Y := \int_0^\infty (|\psi'_z(x)| + |\psi_z(x)| (1 + x^{q+1}) ) \, dx \le C,
$$
uniformly in $z \in B(0,R)$, from what (iv) immediately follows.
\qed

\smallskip
\noindent{\sl Proof of Theorem~\ref{th:CWR}. } Collecting the information obtained on $\LL$ in Corollary~\ref{cor:FirstEV}, Lemma~\ref{lem:Renewal3} \ and using Theorem~\ref{theo:GalWeylSG}, we immediately deduce that \eqref{eq:CWR} holds. \qed

\subsection{Proof of Theorem~\ref{th:MR} in the case without delay}
We present the proof of our main result Theorem~\ref{th:MR} in the case without delay.

\proof[Proof of Theorem~\ref{th:MR} in the case without delay] We split the proof into three steps.

\noindent
{\sl Step 1. A new formulation.} From Lemma \ref{lem:varphig}, for a given  initial datum $0 \le f_0 \in X$ with total densitu of neuron $1$, we may write the solution $f \in C([0,\infty); X)$ to the evolution equation \eqref{eq:ASM} and the solution $F_\lambda$ to the stationary problem \eqref{eq:StSt} as
\bean
\partial_t f + \partial_x f +k(\lambda \varphi[f]) f=0, &\quad& f(t,0) =  \varphi[f(t,\cdot)],\\
\partial_x F +k(\lambda M) F=0, &\quad&  F (0) =  M = \varphi[F],
\eean
 where here and below the $\lambda$ and $x$ dependency is often removed without any confusion.

\smallskip
Next, we consider the variation function $g:= f - F$ which satisfies
\bean
\partial_t g =  -\partial_x g - k(\lambda M) g + (k_\lambda(m) - k_\lambda(M)) f
\eean
complemented with the boundary condition
\bean
g(t,0) &=&  \varphi[f(t,\cdot)] - \varphi[F]\\
 &=& \int_0^\infty k(\lambda \varphi[f]) f  - \int_0^\infty k(\lambda \varphi[F]) F\\
&=&  \MM[g]  +   (\PP_\lambda[f,m] - \PP_\lambda[f,M])  ,
\eean
with $\MM = \MM_\lambda$ defined in \eqref{eq:deMMeps},
Considering the boundary condition as a source term again, we deduce that the variation function $g$ satisfies the equation
\beqn\label{eq:dtg=Lambda+Z}
\partial_t g = \LL g + Z[g],
\eeqn
with the nonlinear term $Z[g] :=  - Q[g] + \delta_0 \QQ[g]$ and where
\bean
Q[g] &:=&  (k_\lambda(\varphi[f]) - k_\lambda(M)) f
\\
\QQ[g] &:=&\PP_\lambda[f,\varphi[f]] - \PP_\lambda[f,M].
\eean

\smallskip
\noindent{\sl Step 2. The nonlinear term.} Using the properties \eqref{eq:th:L1infty}, the assumption \eqref{hyp:a3}{-\eqref{hyp:a4}} and Lemma~\ref{lem:varphig}, we estimate
\bean
\|Q[g]\|_X &=& \|k(\lambda\varphi[f])f-k(\lambda\varphi[F])f \|_{L^1}
\\
&\le& \| f \|_{L^\infty} \, {\xi} \, | \varphi[f] - \varphi[F] | \\
&\lesssim& {\xi} \, \| g \|_{L^1} .
\eean
Similarly, for the boundary term, we have
\bean
|\PP_\lambda [f,m] - \PP_\lambda[f,M]|  \lesssim {\xi} \, \| g \|_{L^1} .
\eean

\smallskip
\noindent{\sl Step 3. Decay estimate.} Thanks to the Duhamel formula, the solution $g$ to the evolution equation \eqref{eq:dtg=Lambda+Z}
satisfies
$$
g(t) = S_{\LL}(t) g_0 + \int_0^t S_{\LL}(t-s)  Z[g(s)] \, \mathrm{d}s.
$$
From Theorem~\ref{th:CWR} and the second step, we deduce
\bean
\| g(t) \|_X
 &\le&  C \, e^{\alpha t} \, \| g_0 \|_X + \int_0^t C \, e^{\alpha (t-s) } \, \| Z[g(s)] \|_X  \, \mathrm{d}s
\\
 &\lesssim&  e^{\alpha t} \, \| g_0 \|_X +  {\xi}\int_0^t  e^{\alpha (t-s) } \, \| g(s) \|_X  \, \mathrm{d}s,
 \eean
for any $t \ge 0$ and for some constant $\alpha < 0$. Thanks to the Gronwall's lemma, we have
\bean
\|g(t)\|_X &\lesssim& e^{\alpha t}\|g_0\|_X+ {\xi}\|g_0\|_X\int_0^t
e^{\alpha t}\exp\{\int_s^t e^{\alpha(t-r)}\mathrm{d}r\}\mathrm{d}s\\
&\lesssim& e^{\alpha t}\|g_0\|_X+{\xi}\, t\,e^{\alpha t}\|g_0\|_X\\
&\lesssim& e^{\alpha' t}\|g_0\|_X,
\eean
for some constant $\alpha' \in (\alpha,0)$. That concludes the proof of Theorem~\ref{th:MR} in the case without delay.
\endproof

\section{Case with delay}
\label{sec:WithDelay}

This section is dedicated to the proof of  Theorem~\ref{th:MR} in the case with delay. Because the arguments are very similar to those
of the previous section, they are only briefly explained.
Following \cite{MSWQ}, we get ride of the delay formulation by writing the problem as a system of PDEs.
We recall or introduce the notations
$$
\PP[h,\mu] := \int_0^\infty k_\lambda(\mu) \, h \, dx, \quad k_\lambda (\mu) = k(x,\lambda\mu),
\quad
\DD[w] := \int_0^\infty w(y) \, b(\mathrm{d}y).
$$
The evolution equation \eqref{eq:ASM} writes
$$
\partial_t f= - \partial_x f - k_\lambda(m)f + \delta_0 p,
$$
with
$$
p (t) = \PP[f(t),m(t)], \quad m(t) = \int_0^\infty p(t-y) \, b(\mathrm{d}y).
$$
Introducing the auxiliary unknown $u$ and the auxiliary equation
$$
\partial_t u = - \partial_y u + \delta_0 p,
$$
the value of $m$ is given by $m(t) = \DD[u(t)]$.
All together, the evolution equation rewrites
\bean
&& \partial_t f= - \partial_x f - k_\lambda(\DD[u])f + \delta_0 p
\\
&& \partial_t u = - \partial_y u + \delta_0 p,
\eean
with $p = \PP[f,\DD[u]]$.

\smallskip
On the other hand, the stationary equation \eqref{eq:StSt} writes
$$
0= - \partial_x F - k_\lambda(M)F+ \delta_0 M
$$
with $M = \PP[F,M]$.
Introducing the auxiliary unknown $U$ and equation
$$
0  = - \partial_y U + \delta_0 M,
$$
the value of $M$ is given by $M = \DD[U]$.
All together, the stationary equation rewrites
\bean
&& 0 =  - \partial_x F - k_\lambda(M)F+ \delta_0 M
\\
&& 0  = - \partial_y U + \delta_0 M,
\eean
with $M = \PP[F,M]$, $M = \DD[U]$.

\smallskip
We define now $g := f - F$, $v = u - U$ and we compute
\bean
\partial_t g &=&  - \partial_x g - k_\lambda(M) g + \bigl( k_\lambda(M) - k_\lambda(m) \bigr) f + \delta_0 \bigl( p(t) - M \bigr)
\\
\partial_tv &=&  - \partial_x v  + \delta_0 \bigl( p(t) - M \bigr)
\eean
with
$$
p(t) - M = \PP[g,M] + \int_0^\infty f \, [k_\lambda(m) - k_\lambda(M)] \, dx.
$$

\smallskip
We first consider the linear system of equations
\bean
\partial_t g &=&  - \partial_x g - k_\lambda(M) g  + \delta_0  \PP[g,M]
\\
\partial_tv &=&  - \partial_x v  + \delta_0 \PP[g,M]
\eean
and the associated semigroup $S_\Lambda$, where $\Lambda$ stands for the operator
$$
\Lambda (g,v) := \bigl( - \partial_x g - k_\lambda(M) g  + \delta_0  \PP[g,M],   - \partial_x v  + \delta_0 \PP[g,M] \bigr).
$$

We introduce the space $X := L^1(\nu) \times L^1(\mu)$ with $\mu := e^{-\delta x} \, dx$ and $\delta>0$ is defined in   condition \eqref{hyp:del}.

\begin{lem}\label{lem:SLambda} The semigroup $S_\Lambda$ associated to $\Lambda$ satisfies
$$
S_\Lambda (t) = \OO(e^{a't}), \quad a '< 0,
$$
in $X$.
\end{lem}

\noindent{\sl Proof of Lemma~\ref{lem:SLambda}. } We write $\Lambda = (\Lambda_1,\Lambda_2)$ and we observe that   $\Lambda_1= \LL$,
where $\LL$ has been defined in the previous section. Because of Theorem~\ref{th:CWR}, we already know that
$$
\| S_{\Lambda_1}(t) (g_0,v_0)\|_{ L^1_q} \lesssim e^{at} \| g_0 \|_{L^1_q}, \quad \forall \, t \ge 0,
$$
for any $(g_0,v_0) \in X$, with $\langle g_0 \rangle = 0$.
Next, we denote $g(t) := S_{\lambda_{_\infty}}(t) (g_0,v_0)$, $v(t) := S_{\Lambda_2}(t) (g_0,v_0)$,   and we compute
\bean
{d \over dt} \int |v(t)| e^{-\delta x}
&=& - \delta \int |v(t)| e^{-\delta x} + | \PP[g(t),M] |
\\
&\le&  - \delta \int |v(t)| e^{-\delta x} + C \, e^{at} \| g_0 \|_{L^1_q}.
\eean
Setting $a' := \max(-\delta,a)$, we conclude by integrating the above differential inequality.
\qed

\smallskip
The last step consists in proving the nonlinear stability by generalizing again the estimates established in the previous section and in \cite{PK1,MSWQ}.
 We write the nonlinear equation as
$$
\partial_t (g,v) =  \Lambda (g,v) + \ZZ
$$
for a nonlinearity/source term $\ZZ =(\ZZ_{_1},\ZZ_2)$ with
$$
\ZZ_{_1}  =  \bigl( k_\lambda(M) - k_\lambda(m) \bigr) f + \delta_0  \int_0^\infty f \, [k_\lambda(m) - k_\lambda(M)] \, dx
$$
and
$$
\ZZ_2  =  \delta_0  \int_0^\infty f \, [k_\lambda(m) - k_\lambda(M)] \, dx.
$$ 
From the assumption \eqref{hyp:a3}-\eqref{hyp:a4}, we have
$$
|\ZZ_2  | \le  \delta_0 \| f \|_{_\infty}    \int_0^\infty|k(x,\lambda m) - k(x,\lambda M) | \, dx \le \delta_0 \, { \xi}\,|m-M|
$$
and similarly
 $$
\|\ZZ_{_1}\|_{L^1} \le { \xi}\,  |m-M|.
$$
From the definition of $\DD$, we have
$$
| m - M | = |\DD[u] - \DD[U]| \le C \, \| v \|_{L^1 (\mu)}.
$$
We finally use Duhamel formula
$$
(g,v) (t) = S_\Lambda(t)(g_0,v_0) +  (S_\Lambda * \ZZ)(t)
$$
and then obtain
$$
\| (g,v)(t) \| \le C_0 e^{at} + C_0\int_0^t e^{a(t-s)} \, { \xi} \, \| (g,v) (s) \| \, ds.
$$
Thanks to the Gronwall lemma, we have
$$
\varphi(t) \le C_0 e^{at} + C_0 \, { \xi} \, \int_0^t e^{a(t-s)} \varphi(s) \, ds =: \psi(t),
$$
which implies
$$
\psi'(t) = a \psi(t) + C_0 \, { \xi} \, \varphi(t) \le (a+ C_0\, { \xi}) \, \psi.
$$
We deduce
$$
\psi(t) \le \psi(0) \, e^{(a+C_0\,{ \xi})t},
$$
and then finally obtain
$$
\varphi(t) \le C_0 \,  e^{(a+C_0\,{ \xi})t}.
$$
We conclude by taking $\lambda > 0$ small and large enough.

\signsm  \signcq  \signqw

\end{document}